\documentclass[11pt,reqno]{amsart}
\usepackage[foot]{amsaddr}
\let\savethebibliography=\thebibliography
\usepackage[square,numbers,comma,sort&compress]{natbib}
\let\thebibliography=\savethebibliography

\usepackage{amsfonts}
\usepackage{url}
\usepackage{color}
\usepackage{amssymb}
\usepackage{enumerate}
\usepackage{fullpage}

\usepackage{verbatim}

\newtheorem{theorem}{Theorem}[section]
\newtheorem{lemma}[theorem]{Lemma}
\newtheorem{Lemma*}[theorem]{Lemma}

\newtheorem{proposition}[theorem]{Proposition}
\newtheorem{corollary}[theorem]{Corollary}

\theoremstyle{definition}

\newtheorem{example}[theorem]{Example}

\newcommand{\Q}{{\mathbb Q}}
\newcommand{\Z}{{\mathbb Z}}
\newcommand{\N}{{\mathbb N}}

\newcommand{\R}{\mathbb R}

\DeclareMathOperator\st{s.t.}
\DeclareMathOperator\argmin{argmin}

\DeclareMathOperator\BMILP{BMILP}
\DeclareMathOperator\BILP{BILP}

\DeclareMathOperator\proj{proj}

\DeclareMathOperator\cl{cl}

\usepackage{ifthen}
\makeatletter
\newcommand{\DeclareBracket}[3]{
  \newcommand{#1}[2][]{%
  \ifthenelse%
  {\equal{##1}{}}%
  {\left#2##2\right#3}%
  {\csname ##1l\endcsname#2##2\csname ##1r\endcsname#3}}}
\makeatother
\DeclareBracket\abs||           
\DeclareBracket\norm\|\|        
\DeclareBracket\floor\lfloor\rfloor
\DeclareBracket\ceil\lceil\rceil
\DeclareBracket\set\{\}
\DeclareBracket\paren()
\DeclareBracket\bracket[]
\DeclareBracket\inner\langle\rangle
\DeclareBracket\fractional\{\}
\makeatother

\usepackage{amssymb}

\usepackage{enumerate}
\def\ve#1{\mathchoice{\mbox{\boldmath$\displaystyle\bf#1$}}
{\mbox{\boldmath$\textstyle\bf#1$}}
{\mbox{\boldmath$\scriptstyle\bf#1$}}
{\mbox{\boldmath$\scriptscriptstyle\bf#1$}}}

\usepackage{ifpdf}

\usepackage[breaklinks=true,colorlinks,citecolor=blue,linkcolor=blue]{hyperref}

\numberwithin{equation}{section}

\title[Parametric integer programming algorithm for bilevel mixed integer programs]{A
    parametric integer programming algorithm for bilevel mixed integer programs}

\author[M.~K\"oppe]{Matthias K\"oppe}
\address{M.~K\"oppe: University of California, Davis, Department of
  Mathematics, One Shields Avenue, Davis, CA 95616, USA}
\email{mkoeppe@math.ucdavis.edu}

\author[M.~Queyranne]{Maurice Queyranne}
\address{M.~Queyranne:
Centro de Modelamiento Matem\'atico, Unit\'e Mixte Internationale du CNRS (France),
Universidad de Chile, Av.~Blanco Encalada 2120, Santiago, Chile;
and Sauder School of Business at the University of British Columbia,
2053 Main Mall, Vancouver, BC, Canada, V6T 1Z2}
\email{maurice.queyranne@sauder.ubc.ca}
\thanks{The research of the last two authors was supported in part by a Discovery grant from the Natural Sciences and Engineering Research Council (NSERC) of Canada to the second author.}

\author[C.~T.~Ryan]{Christopher Thomas Ryan}
\address{C.~T.~Ryan: Sauder School of Business at the University of British Columbia,
2053 Main Mall, Vancouver, BC, Canada, V6T 1Z2}
\email{chris.ryan@sauder.ubc.ca}

\date{\today}

\begin{document}

\maketitle

\subsection*{Abstract}
We consider discrete bilevel optimization problems where
the follower solves an integer program with a fixed number of variables.
Using recent results in parametric integer programming,
we present polynomial time algorithms for pure and mixed integer bilevel problems.
For the mixed integer case where the leader's variables are continuous,
our algorithm also detects whether the infimum cost fails to be attained,
a difficulty that has been identified but not directly addressed in the literature.
In this case it yields a ``better than fully polynomial time'' approximation scheme
with running time polynomial in the logarithm of the absolute precision.
For the pure integer case where the leader's variables are integer,
and hence optimal solutions are guaranteed to exist,
we present an algorithm which runs in polynomial time when the total number of variables is fixed.

\subsection*{Keywords} bilevel mixed integer linear programming, parametric integer linear programming, computational complexity, binary search

\section{Introduction}\label{s:intro-fixed}
Bilevel integer programs are two-stage decision problems in which two decision
makers act in sequence and where some or all of the decision variables are integer.
A \emph{leader} first decides on some variables, denoted here by $\ve z$
(called the leader's, or upper-level, variables),
which influence the decision of a \emph{follower},
who then chooses $\ve x$ (the follower's, or lower-level, variables).
The leader is interested in optimizing over the joint decision $(\ve x, \ve z)$
subject to joint resource and incentive compatibility constraints.
This setting has been used to model a large variety of applied problems,
including toll setting in transportation networks \cite{Labbe1998},
revenue management \cite{Cote2003},
competitive location planning \cite{Fischer2004-facility-bilevel}, and
vehicle routing \cite{marinakis2007new} among many others.

To date the majority of literature has focused on the setting where
all constraints and objectives are linear and
both the leader and follower's variables are continuous (this is evident in observing,
for instance, the survey papers \cite{chinchuluun2009multilevel, Colson2007}, research monograph \cite{DempeText} and
journal special issue \cite{migdalas1996}). Significant work has also been undertaken
in nonlinear continuous settings (see for instance \cite{al1992global,migdalas1995nonlinear}).
Settings where some or all variables are integer are less well-studied, but
have nonetheless garnered interest in recent years \cite{Dempe2001,Dempe2000-knapsack,Denegre2008,Gumus2005,Vicente1996}.

The main focus of this paper is the following
bilevel mixed integer linear program (BMILP),
where the leader's variables $\ve z\in\R^d$ are continuous
and the follower's variables $\ve x\in\Z^n$ are integer:
\begin{eqnarray}
    \inf_{\ve x, \ve z} && \ve c \cdot \ve x + \ve e \cdot \ve z \label{eq:general-fixed_obj}\\
    \st && C\ve x + D\ve z \le \ve p,\ \ \ve z \ge \ve 0 \label{eq:general-fixed_P}\\
    && \ve x \in \argmin_{\ve x'}\{\, \ve \psi \cdot \ve x':
            A \ve x' \leq B\ve z + \ve u,\ \ \ve x'\in\Z^n \,\}\label{eq:general-fixed_argmin}
\end{eqnarray}
\noindent
where all data are integer.

An interpretation of the problem is that the leader chooses $\ve z$
and as a result affects
the right-hand sides of the lower level (follower's) problem
\[
  \min_{\ve x'} \{\ve \psi\cdot\ve x' : A\ve x'\le B\ve z + \ve u,\ \ \ve x'\in\Z^n\},
\]
which is a \emph{parametric integer program} with right-hand side parameterized by~$\ve z$.
This interpretation underscores the parametric nature of bilevel integer programs that
will be a major focus in this paper.

An important comment on the formulation is the implicit assertion of what is known
in the bilevel literature as the ``optimistic assumption'' (see for instance \cite{Colson2007}).
It is possible for the lower level problem to have multiple optimal
solutions for a given $\ve z$.
The follower is indifferent amongst these
alternate optima, but the leader may prefer certain optima over others, and yet can exert
no more direct influence over the follower's choice.
How should the follower decide which
of the alternate optima to choose?
To eliminate this difficulty we make the \emph{optimistic assumption}:
whenever the follower faces alternate optima for a given $\ve z$,
the follower chooses an optimum $\ve x$ that best suits the leader.
Indeed, this is reflected in
the formulation \eqref{eq:general-fixed_obj}--\eqref{eq:general-fixed_argmin} itself:
the leader is modeled as choosing $\ve x$ within the $\argmin$ set of the lower level problem.

An important object studied throughout is the set of
\emph{bilevel feasible solutions}
\begin{equation}\label{eqn:bilevel_feasible_optimistic}
  \mathcal{F} = \left\{ (\ve x, \ve z) \in \R^{n+d} :
        \eqref{eq:general-fixed_P},\ \eqref{eq:general-fixed_argmin}
        \right\}.
\end{equation}
We assume that this \emph{bilevel feasible set} $\mathcal{F}$ is bounded (though possibly empty). Much of the complication associated with bilevel integer programming is captured in
the properties of this set. Firstly, $\mathcal{F}$ is known to be nonconvex in general,
and difficult to describe directly \cite{Colson2007}. Secondly, given $\ve z$ is continuous
it is possible that $\mathcal{F}$ is not closed \cite{Vicente1996}.
This possibility is reflected in our definition of BMILP in \eqref{eq:general-fixed_obj}--\eqref{eq:general-fixed_argmin} where the leader's optimization is expressed as an infimum and not a minimum. When $\ve z$ is also integer and the instance is feasible, then an optimal
solution always exists and ``$\inf$'' can be replaced with ``$\min$'' \cite{Vicente1996}.

\begin{example}
To demonstrate this complication, consider the following simple example.
\begin{eqnarray*}
    \inf_{x, z} && -x+z\\
    \st && 0 \le z \le 1\\
    && x \in \argmin_{x'}\{\, x':
            x' \ge z,\  0 \le x' \le 1, \ \ x'\in\Z \,\}
\end{eqnarray*}
The ``$\argmin$" then simply requires that $x = \lceil z \rceil $ and thus the problem reduces to:
\begin{equation*}
    \inf_{z} \{ z - \lceil z \rceil : 0 \le z \le 1\}.
\end{equation*}
Clearly, the infimum is $-1$ and cannot be attained. Observe that the set of bilevel feasible solutions
$$\mathcal{F} = \{(0,0)\} \bigcup (\{1\} \times (0,1])$$
is neither convex nor closed.
\end{example}

The approaches to problem \eqref{eq:general-fixed_obj}--\eqref{eq:general-fixed_argmin}
in the literature use traditional ideas
from integer optimization, such as branch-and-bound enumeration
\cite{MooreBard1990} or cutting planes \cite{Denegre2008}, all based on continuous relaxations. It is well
documented that there are numerous pitfalls when applying standard ideas of branch-and-bound
or cutting plane methods to a bilevel setting (see \cite{Denegre2008} for a discussion).
Moreover, the algorithm introduced by Moore and Bard in \cite{MooreBard1990} applies
only to the case that both $\ve x$ and $\ve z$ are integer,
or otherwise where it is known that an optimum solution exists.
Thus, in general, this method
cannot solve \eqref{eq:general-fixed_obj}--\eqref{eq:general-fixed_argmin},
in the sense that it cannot distinguish when
an infimum is attained or not.
The algorithm of DeNegre and Ralphs \cite{Denegre2008} deals exclusively with
the pure integer setting.
In addition, these two papers give no run-time guarantees.

Our paper makes the following contributions:
\begin{enumerate}[(i)]
\item We present (Theorem~\ref{thm:main-result-ES}) an algorithm which solves BMILP
\eqref{eq:general-fixed_obj}--\eqref{eq:general-fixed_argmin}
in full generality (under our boundedness and integral data assumptions)
by identifying when the infimum is attained and finding
an optimal solution if one exists.
No previously proposed algorithm is capable of this.
\item Our algorithms runs in polynomial time when
the follower's dimension $n$ is fixed for the mixed integer case (i.e., $z$ is continuous),
and when the total dimension $n+d$ is fixed in the pure integer case (i.e., when $z$ is integer).
To the authors' knowledge, these are the
first polynomial run-time guarantees of this type for discrete bilevel optimization. These results
can be seen as a common extension of two important results on optimization in fixed dimension. The
first is due to Lenstra \cite{Lenstra83}, who establishes a polynomial time algorithm for integer
programming with a fixed
number of variables. The second is due to Deng \cite{Deng1997-complexity}, who develops a polynomial
time algorithm for bilevel \emph{linear} programs when the number of decision variables of the
follower is fixed. Our result extends Lenstra's result to a bilevel setting and Deng's result to
an integer setting. For further discussion of complexity issues in bilevel optimization see
\cite{chinchuluun2009multilevel}.
\item When no optimal solution to BMILP exists, we find,
for any given $\epsilon > 0$ and in time polynomial in the input size and $\log(1/\epsilon)$,
a solution with objective value no more than $1+\epsilon$ times the infimum.
\end{enumerate}

In Section~\ref{s:def}, we introduce formal definitions and prove our first key result:
when an instance of BMILP is feasible
its infimum value is a rational number with polynomially bounded binary encoding size.
In Section~\ref{s:parametric}, we apply recent results from the study of parametric
integer programming by Eisenbrand and Shmonin~\cite{Eisenbrand2008}, based on the work of Kannan~\cite{Kannan1992-lattice}.
These references provide polynomial time algorithms
for decision versions of parametric integer programs.
We apply these algorithms in the mixed integer bilevel setting and,
using binary search, derive our main result (Theorem~\ref{thm:main-result-ES}):
a polynomial time algorithm for the mixed integer case where $\ve z$ is continuous and
the follower's dimension $n$ is fixed.
Finally, in Section~\ref{s:pure-integer} we present two algorithms to find optimal solutions
to the pure integer case; that is, BMILP under the further restriction that $\ve z$ is integer. Both of these algorithms run
in polynomial time when the total dimension $n+d$ is fixed. The first algorithm is similar
to those discussed in the previous section, whereas the second is based on an approach similar to that found in  \cite{Koppe2008-games}
which is based on a theory of rational generating functions.

\section{Definition and Preliminary Results}\label{s:def}
An instance of the bilevel mixed integer linear program (BMILP)
\eqref{eq:general-fixed_obj}--\eqref{eq:general-fixed_argmin}
is defined by the leader's and follower's dimensions $d$ and $n$;
the numbers $h$ and $m$ of linear inequality constraints
in~\eqref{eq:general-fixed_P}
and in the feasible set of the follower's subproblem~\eqref{eq:general-fixed_argmin};
integer matrices $A \in \Z^{m\times n}$,
$B \in \Z^{m\times d}$, $C \in \Z^{h\times n}$ and $D \in \Z^{h\times d}$;
integer objective coefficient vectors $\ve c \in \Z^n$, $\ve e \in \Z^d$
and $\ve\psi \in \Z^n$;
and integer right-hand side vectors
$\ve u \in \Z^m$ and $\ve p \in \Z^h$.
The binary encoding input size of an instance is
the total number of bits needed to input all the data defining the instance.

We assume that the upper level polyhedron
\begin{equation}\label{eq:def-P}
     P = \{(\ve x,\ve z)\in\R^{n + d} : \eqref{eq:general-fixed_P}\}
\end{equation}
is bounded (but possibly empty).
We also assume that for every
$\ve z\in\proj_{\ve z} P$, where
\[
  \proj_{\ve z} P = \{\ve z\in\R^d : (\ve x,\ve z)\in P \text{ for some }\ve x\in\R^n\},
\]
the lower level feasible set $\{(\ve x,\ve z)\in\R^{n + d} : A\ve x \le B\ve z +\ve u\}$
is also bounded (and possibly empty).
(For example, by Farkas's Lemma,
if lower bounds $\ve x \ge \underline{\ve x}$
and $\ve z \ge \underline{\ve z}$
are part of or implied by the constraints defining~$P$,
and also by those defining the lower level feasible set, then
the boundedness assumptions are satisfied
if and only if the $(n+d)$-row vector of all 1's
is a nonnegative combination of the rows of the constraint matrix $(C,D)$,
and also a nonnegative combination of the rows of the matrix $(A,-B)$.)

Under these boundedness assumptions,
an instance of BMILP is either infeasible or else must have a bounded infimum value.
Feasibility can be decided in polynomial time
by applying Lenstra's (mixed) integer programming algorithm \cite{Lenstra83}
to decide if
$\{(\ve x,\ve z)\in P : A \ve x \leq B\ve z + \ve u,\ \ \ve x\in\Z^n\} \not=\emptyset$.
A key result used in the sequel is that, when the instance of BMILP is feasible,
its infimum value is a rational number
with polynomially-bounded binary encoding size.

\begin{proposition}\label{prop:denominator-bound}
  Assume that
  a given instance of
  problem \eqref{eq:general-fixed_obj}--\eqref{eq:general-fixed_argmin} is feasible and
  has a bounded optimum (infimum) value.
  Then the infimum is a rational number with denominator bounded above by
  the maximum absolute value of a sub-determinant of the matrix
  $(D,B)^T$.
\end{proposition}
\begin{proof}
Let $(\bar {\ve x},\bar {\ve z})\in\mathcal{F}$ be a bilevel feasible solution
to \eqref{eq:general-fixed_obj}--\eqref{eq:general-fixed_argmin} and define
\[
  Q(\bar {\ve x},\bar {\ve z}) = \{\ve z\in\R^d :  D\ve z \le \ve p - C\bar {\ve x};
    \ \ \beta_i(\bar {\ve z}) \le B_i\ve z + u_i < \beta_i(\bar {\ve z}) +1\ \forall i\in M\},
\]
where $M=\{1,\dots,m\}$, $B_i$ is the $i$-th row of~$B$, and
$\beta_i(\bar {\ve z}) = \lfloor B_i \bar {\ve z} + u_i \rfloor$.
Note that $\bar {\ve z}\in Q(\bar {\ve x},\bar {\ve z})
$,
so $Q(\bar {\ve x},\bar {\ve z})$ is a nonempty, bounded quasi-polyhedral set,
i.e., its topological closure is
\begin{equation}\label{eq:def-clQxz}
 \cl Q(\bar {\ve x},\bar {\ve z}) = \{\ve z\in\R^{d} :  D\ve z \le \ve p - C\bar{\ve x};
   \ \ \beta_i(\bar {\ve z}) \le B_i \ve z + u_i \le \beta_i(\bar {\ve z}) +1\ \forall i\in M\},
\end{equation}
a nonempty polytope.
Every $\ve z\in Q(\bar {\ve x},\bar {\ve z})$ satisfies
$\lfloor B_i \ve z + u_i \rfloor = \lfloor B_i \bar {\ve z} + u_i \rfloor$ for all $i=1,\dots,m$ and,
since $A$ is integer,
\[
 \forall \ve x\in \Z^n: \quad A\ve x \le B\ve z + \ve u \Longleftrightarrow A\ve x \le
 B\bar {\ve  z} + \ve u,
\]
implying that
$\bar {\ve x} \in \argmin\{\ve \psi \cdot \ve x : A\ve x \le B\ve z + \ve u,\ \ve x \in \Z^n\}$
for all $\ve z\in Q(\bar {\ve x},\bar {\ve z})$,
that is,
\begin{equation}\label{eq:xbar}
  \forall \ve z\in Q(\bar {\ve x},\bar {\ve z})\quad (\bar {\ve x},\ve z)\in\mathcal{F}.
\end{equation}
Now let $v^*$ denote the infimum value in BMILP \eqref{eq:general-fixed_obj}--\eqref{eq:general-fixed_argmin}.
There is a sequence $(\hat{\ve x}^k, \hat{\ve z}^k)_{k\in\N}$
of bilevel feasible solutions with objective value
$\ve c \cdot \hat{\ve x}^k + \ve e \cdot \hat{\ve z}^k$ converging to~$v^*$.
Since $P$ is a polytope, its projection
$\proj_{\ve x}P = \{\ve x\in\R^n : (\ve x,\ve z)\in P \text{ for some }\ve z\in\R^d\}$
contains a finite number of integer points $\ve x\in P\cap\Z^n$.
Therefore there is a subsequence $(\tilde {\ve x}^k, \tilde {\ve z}^k)_{k\in\N}$
with all $\tilde{\ve x}^k = \tilde{\ve x}^0$.
For every $i\in M$ there is also a finite number of possible integer values
$\beta_i(\ve z) = \lfloor B_i \ve z + u_i \rfloor$
when $\ve z\in\proj_{\ve z} P$.
%
Thus there exists a subsequence $(\tilde{\ve x}^0, \breve{\ve z}^k)_{k\in\N}$ in $\mathcal{F}$
with all $\beta(\breve{\ve z}^k) = \beta(\breve{\ve z}^0)$,
i.e., $(\tilde{\ve x}^0, \breve{\ve z}^k)\in Q(\tilde{\ve x}^0, \breve{\ve z}^0)$,
and $\ve c\cdot\tilde{\ve x}^0 + \ve e\cdot\breve{\ve z}^k$ converging to~$v^*$.
From (\ref{eq:xbar}) it follows that
$v^* = \ve c \cdot \tilde{\ve x}^0 +
  \inf_{\ve z}\{\ve e \cdot \ve z : \ve z\in Q(\ve {\tilde x}^0,\breve{\ve z}^0) \}$.
By continuity of the objective~$\ve e \cdot \ve z$,
$v^* - \ve c\cdot\tilde{\ve x}^0$ is the optimum value of the linear programming problem
$\min_{\ve z}\{\ve e\cdot \ve z : \ve z\in \cl Q(\tilde{\ve x}^0,\breve{\ve z}^0)\}$.
By~(\ref{eq:def-clQxz}),
for every $i\in M$ at most one of the two constraints
$\beta_i(\breve{\ve z}^0) \le B_i \ve z+ u_i$ or  $B_i \ve z + u_i\le \beta_i(\breve{\ve z}^0) +1$
is binding at the LP optimum, hence $v^* - \ve c\cdot\tilde{\ve x}^0$ must also be the optimum value of
a linear programming problem
\begin{equation}\label{eq:slice-lp}
  \begin{aligned}
    \min\quad & \ve e\cdot\ve z &\\
    \st\quad  & D\ve z \le \ve p - C\tilde{\ve x}^0 &\\
    & -B_i \ve z \le u_i -\beta_i(\breve{\ve z}^0), &  \forall i\in I\\
    & B_j \ve z  \le -u_j + \beta_j(\breve{\ve z}^0) + 1, & \forall j\in J
  \end{aligned}
\end{equation}
for some disjoint $I,\;J\subseteq M$.
It is well known that the optimal value of a bounded linear program with integer data
occurs at a basic feasible solution, the denominators of which are bounded by
the maximum absolute value of a sub-determinant of its constraint matrix.
By simple properties of determinants the stated result then follows.
\end{proof}

Note that in the pure integer case,
where the leader's variables are also required to be integer,
the bilevel feasible set is
\begin{equation}\label{eqn:pure_bilevel_feasible_optimistic}
  \mathcal{F'} = \left\{ (\ve x, \ve z) \in P : \ve z\in\Z^d,\ \ \ve x \in \argmin\{\,
        \ve \psi \cdot \ve x' : A \ve x' \leq B\ve z + \ve u, \,\ve x'\in\Z^n \,\} \right\}.
\end{equation}
When it is non-empty, $\mathcal{F'}$ is a finite set, and then
the infimum in \eqref{eq:general-fixed_obj} is attained.
Since $\ve c$ and $\ve e$ are integral, this optimum value is integral,
providing a trivial strengthening of Proposition~\ref{prop:denominator-bound} in the pure integer case.

\section{An Algorithm for Bilevel Mixed Integer Linear Programming}\label{s:parametric}
In this section we prove the main result of this paper:
\begin{theorem}\label{thm:main-result-ES}
There exists an algorithm which
solves the BMILP problem \eqref{eq:general-fixed_obj}--\eqref{eq:general-fixed_argmin}
in the following sense: the algorithm
\begin{enumerate}[(i)]
\item decides if the instance is feasible;
\item if it is feasible, decides if an optimal solution exists (that is, if the infimum is attained); and
\item if the infimum is attained, finds an optimal solution.
\end{enumerate}
The algorithm runs in polynomial time when the follower's dimension $n$ is fixed.
\end{theorem}

To prove this theorem, we use results on parametric integer programs
with parameterized right-hand sides.
Parametric integer programming is well-studied,
the particular approach we discuss builds on geometric and algorithmic number theory and
has it roots in the work of Lenstra \cite{Lenstra83} on integer linear programming in
fixed dimension.
We draw mainly on the recent work of Eisenbrand and Shmonin \cite{Eisenbrand2008}, which
itself builds on earlier results by Kannan \cite{Kannan1992-lattice}.
These two related works solve a decision version of parametric integer programming (defined below) in polynomial time when certain parameters are fixed.
As discussed in Section~\ref{s:intro-fixed}, we can view the lower level problem
as a parametric integer program.
But the connection between parametric integer programming and bilevel integer programming
in fact goes deeper, as we now discuss.

The \emph{integer projection} of a polyhedron $Q\subseteq\R^{m + p}$ is the set
\begin{equation}\label{eq:def-integer-proj}
   Q/\Z^p = \{\ve b\in\R^m : (\ve b,\ve w)\in Q \text{ for some }\ve w\in\Z^p\}.
\end{equation}
The feasibility version of parametric integer linear programming is
the following decision problem:
\begin{quotation}
    \textsc{PILP}: Given a rational matrix $R \in \Q^{m \times n}$ and
    a rational polyhedron $Q\in \R^{m + p}$ decide if
    the linear system $R\ve y \le \ve b$ has an integer solution for all $\ve b\in Q/\Z^p$.
\end{quotation}
Letting
\[
   P_{\ve b} = \{\ve y \in \R^n: R\ve y \le \ve b\},
\]
this problem is equivalent, by negation, to deciding the sentence
\begin{equation}\label{eq:PILP}
\exists \ve b \in Q / \Z^p \ \ P_{\ve b}\cap\Z^n = \emptyset.
\end{equation}
Eisenbrand and Shmonin \cite{Eisenbrand2008} prove:
\begin{theorem} \label{thm:PILP-algorithm} \emph{(Theorem 4.2 in \cite{Eisenbrand2008})}
There is an algorithm that, given a rational matrix
$R \in \Q^{m \times n}$ and a rational polyhedron $Q \subset \R^{m+p}$, decides the sentence
\eqref{eq:PILP}.
The algorithm runs in polynomial time if $p$ and $n$ are fixed.
\end{theorem}
%
%
%
%
%

%
%
Our goal is to establish a connection between \textsc{PILP} and BMILP so this result can
be used in solving BMILP.
First we define a decision version  as follows:
\begin{quotation}\label{prob:BMILP-alpha}
   $\BMILP_{\alpha}$: Given integer matrices $A \in \Z^{m\times n}$,
   $B \in \Z^{m\times d}$, $C \in \Z^{h\times n}$ and $D \in \Z^{h\times d}$;
   integer vectors $\ve c \in \Z^n$, $\ve e \in \Z^d$, $\ve\psi \in \Z^n$,
   $\ve u \in \Z^m$ and $\ve p \in \Z^h$;
   and rational scalar $\alpha\in\Q$,
   decide if there exists a vector $(\ve x, \ve z) \in \R^{n+d}$ such that
   \eqref{eq:general-fixed_P}, \eqref{eq:general-fixed_argmin} and
   $\ve c \cdot \ve x + \ve e \cdot \ve z \le \alpha$.
\end{quotation}
\begin{proposition}\label{prop:decision-version}
There exists an algorithm that decides $\BMILP_{\alpha}$.
The algorithm runs in polynomial time when the follower's dimension $n$ is fixed.
\end{proposition}
\begin{proof}
Since $\ve \psi$ is integral, we can express \eqref{eq:general-fixed_argmin}
as the fact that
\begin{equation}\label{eq:x-is-argim}
  \{\ve x' \in \R^n :  \ve\psi\cdot\ve x' \le \ve\psi\cdot\ve x - 1,
                   \ \ A\ve x' \le B\ve z + \ve u\}\cap\Z^n
    = \emptyset.
\end{equation}
Let 
$R = (\ve \psi, A, 0_{d\times n})^T$
where $0_{d\times n}$ is a $d\times n$
matrix of zeroes;
$\ve b = (b_{(1)}, \ve b_{(2)}, \ve b_{(3)})^T$
with $b_{(1)} = \ve\psi\cdot\ve x - 1$,
     $\ve b_{(2)} = B\ve z + \ve u$ and
     $\ve b_{(3)} =  \ve z$.
Then \eqref{eq:general-fixed_argmin} and $\ve z\ge \ve 0$ together are
equivalent to $P'_{\ve b}\cap\Z^n = \emptyset$ where
$P'_{\ve b} = \{\ve x' \in \R^n: R\ve x' \le \ve b\}$.
Letting
\begin{eqnarray}
Q_{\alpha} = \{ (\ve b, \ve x) &:&
         \ve c\cdot\ve x + \ve e\cdot\ve b_{(3)}\le \alpha,
     \ \ C\ve x + D\ve b_{(3)} \le\ve p,
     \ \ \ve b_{(3)} \ge \ve 0,             \nonumber\\
     &&  A\ve x \le \ve b_{(2)},
     \ \  b_{(1)} = \ve\psi\cdot\ve x - 1,
     \ \ \ve b_{(2)} = B\ve b_{(3)} + \ve u \label{eq:def-Q-alpha}
  \},
\end{eqnarray}
deciding $\BMILP_{\alpha}$ reduces to deciding an instance of~\eqref{eq:PILP} with
$p = n$, $\ve w = \ve x$, $Q = Q_{\alpha}$
and $P_{\ve b} = P'_{\ve b}$.
The result now follows from Theorem~\ref{thm:PILP-algorithm}.
\end{proof}

Our approach to the proof of Theorem~\ref{thm:main-result-ES}
is to use binary searches on the target value $\alpha$ in $\BMILP_{\alpha}$
and related decision problems, so as
to find the infimum value of problem BMILP and an optimal solution
if one exists.


%
%
%
The following result from
Kwek and Mehlhorn~\cite{Kwek2003-rational-search}
yields the binary search algorithm used in proving our results:
\begin{theorem}\label{thm:kwek-rational-search} \emph{(Theorem 1 in \cite{Kwek2003-rational-search})}
There exists a $\theta(\log L)$-time
algorithm that determines an unknown rational number
$r$ which lies in the set $\mathcal{Q}_L = \left\{ \frac {p}{q}: p,q \in \{1,\ldots, L\}\right\}$
by asking at most $2\log_2 L + O(1)$ queries of the form ``is $r \le \alpha$?''.
\end{theorem}
%
%

We first establish:
\begin{lemma}\label{lemma:infimum-value}
Consider a feasible instance to
BMILP problem \eqref{eq:general-fixed_obj}--\eqref{eq:general-fixed_argmin}.
Its infimum value can be determined in polynomial time when the dimension $n$ is fixed.
\end{lemma}
\begin{proof}
We first determine polynomial-sized lower and upper bounds $\underline{v}$ and $\bar v$ on
the infimum value $v^*$
by solving in polynomial time the linear programs
$\underline{v} = \min\{\ve c\cdot\ve x + \ve e\cdot\ve z : (\ve x,\ve z)\in P\}$
and
$\bar v        = \max\{\ve c\cdot\ve x + \ve e\cdot\ve z : (\ve x,\ve z)\in P\}$.
We then perform ordinary binary search for $v^*$ in successive intervals $[v', v'']$
satisfying $v' < v^* \le v''$,
where the search query ``is $v^* \le \alpha$'' is the decision problem $\BMILP_\alpha$,
and the initial interval $[v', v''] =[ \underline{v}-1, \bar v]$.
We stop, after a polynomial number of queries, as soon as $v'' - v' <1$.
If there exists an integer $\beta$ such that $v' < \beta < v''$ then we perform one more query
``is $v^* \le \beta$'' and set the current interval to $[\beta-1, \beta]$ if the response was positive
(i.e., if $v^* \le \beta$), and to $[\beta, v'']$ otherwise.
If there is no such integer $\beta$ then we set the current interval to $[\lfloor v'\rfloor, v'']$.
In any case we obtain an interval $[\bar v', \bar v'']$ such that
$\bar v' < v^* \le\bar v''$, $\bar v'' - \bar v' \le 1$
and $\bar v'$ is integer.
We now invoke the binary search algorithm of Theorem~\ref{thm:kwek-rational-search}
for the rational number $r = v^* - \bar v'\in\mathcal{Q}_L$ where $L$ is
the upper bound on the denominator of $v^*$
from Proposition~\ref{prop:denominator-bound}.
\end{proof}

We are now ready to prove Theorem~\ref{thm:main-result-ES}.

\begin{proof}[Proof of Theorem~\ref{thm:main-result-ES}]
(i) As indicated in Section~\ref{s:def}, feasibility can be decided in polynomial time
by applying Lenstra's (mixed) integer programming algorithm \cite{Lenstra83}
to decide if
$\{(\ve x,\ve z)\in P : A \ve x \leq B\ve z + \ve u,\ \ \ve x\in\Z^n\} \not=\emptyset$.

(ii) Assume the instance is feasible.
By Lemma~\ref{lemma:infimum-value} we can find the infimum value $v^*$ in polynomial time.
This infimum is attained if and only if there exists a bilevel feasible $(\ve x, \ve z) \in \mathcal{F}$ such that
$\ve c \cdot \ve x + \ve e \cdot \ve z = v^*$.
This can be decided in polynomial time as in the proof of Proposition~\ref{prop:decision-version}
by replacing the first constraint in the definition of $Q_\alpha$
in equation~\eqref{eq:def-Q-alpha}
with $\ve c \cdot \ve x + \ve e \cdot \ve b_{(3)} = v^*$.

(iii) Finally, assume that the infimum value is attained.
We extend to bilevel mixed integer programming a standard method used in
(ordinary) mixed integer programming to successively determine an optimum solution
given an optimum objective value oracle (see for instance \cite{Schrijver1986}, page 260).
In this method, we
first
find the optimum value of related instances of BMILP to successively determine
each component of an optimal solution $(\ve x^*, \ve z^*)$
with lexicographically smallest integer component $\ve x^*$.
By induction,
assume that all components $x_k^*$ for $k = 1,\dots, j-1$ have already been determined.
Then $x^*_j$ is the optimum value of the BMILP
\begin{eqnarray*}
    \min_{\ve x, \ve z} && x_j \\
    \st &&     C\ve x + D\ve z \le \ve p;\ \ \ve z \ge \ve 0;
            \ \ \ve c \cdot \ve x + \ve e \cdot \ve z = v^*;\\
        &&     x_k = x_k^*\;\text{ for }k = 1,\dots,j-1; 
            \ \ \eqref{eq:general-fixed_argmin}.
\end{eqnarray*}
Since this problem is feasible and its objective is integer valued and bounded,
its optimum value is attained, and
is determined in polynomial time by Lemma~\ref{lemma:infimum-value}.
(Note that, since this optimal value $x^*_j$ is integer,
we only need to use ordinary binary search and
stop when the width of the current interval is less than one.)

Having determined $\ve x^*$, we next determine the continuous component $\ve z^*$.
Note that, even though $\ve x^*$ has now been fixed,
$\ve z^*$ cannot be directly found by solving a simple linear program,
because we still need to enforce
the (nonconvex) incentive compatibility constraint~\eqref{eq:general-fixed_argmin}.
%
%
Using ideas from the proof of Proposition~\ref{prop:denominator-bound},
we construct an ``integer RHS vector'' $\ve r = (\lfloor B_i \ve z + u_i\rfloor)_{i=1,\dots m}$
which is lexicographically minimum among all those defined by optimum BMILP solutions $(\ve x^*, \ve z)$.
Thus by induction assume that for $k = 1,\dots, i-1$ we have determined integers $r_k$
such that 
the set
\[ %
R_{i-1} = \left\{(\ve x,\ve z) \in \R^{n+d}:
    \begin{array}{lll}
        C\ve x^*+D\ve z\le\ve p; \ \ \ve z \ge \ve 0;\ \ \ve c \cdot \ve x^*+\ve e\cdot \ve z = v^*; \\
        r_k \le B_k\,\ve z + u_k < r_k + 1\text{ for all }k < i;
            \ \ \ve x = \ve x^*;\ \ \eqref{eq:general-fixed_argmin}
    \end{array}\right\}
\]
is nonempty.
As in the proof of Proposition~\ref{prop:denominator-bound},
note that the closure $\cl R_{i-1}$ is determined by replacing each strict inequality
$B_k\ve z + u_k < r_k + 1$ with $B_k\ve z + u_k \le r_k + 1$,
and is thus a nonempty polytope.
We then compute the infimum value of the BMILP
\begin{equation}\label{eq:BMILP-rho_i}
    \rho_i = \inf_{\ve x, \ve z}\{ B_i\,\ve z + u_i : \ve z\in \cl R_{i-1}\},
\end{equation}
in polynomial time by Lemma~\ref{lemma:infimum-value}.
Let $r_i = \lfloor \rho_i \rfloor$.
We claim that $R_i\not=\emptyset$.
To prove this claim, let $\big((\ve x^*,\ve{\bar z}^j)\in R_{i-1}\big)_{j\in\N}$ be
a sequence of points in~$\cl R_{i-1}$
with $\lim_{j\rightarrow\infty}(B_i\,\ve{\bar z}^j + u_i) = \rho_i$.
By continuity of the linear function $(\ve x^*,\ve z)\mapsto B_i\,\ve z + u_i$,
there exists a sequence $\big((\ve x^*,\ve z^j)\in R_{i-1}\big)_{j\in\N}$
of points in~$R_{i-1}$ with
$\lim_{j\rightarrow\infty}(B_i\,\ve z^j + u_i) = \rho_i < \lfloor\rho_i\rfloor + 1$.
This implies that $\lfloor\rho_i\rfloor \le B_i\,\ve z^j + u_i < \lfloor\rho_i\rfloor + 1$
for some $j\in\N$, and therefore $(\ve x^*,\ve z^j)\in R_i$,
showing that $R_i\not=\emptyset$.

Thus, after determining $\ve x^*$ and
computing the infimum value $\rho_i$ of $m$  BMILPs (\ref{eq:BMILP-rho_i}),
we obtain the polynomial-sized, integer vector $\ve r = \big(\lfloor \rho_i \rfloor\big)_{i=1,\dots,m}$
such that there exists an optimal solution $(\ve x^*, \ve z)$ to the original
BMILP \eqref{eq:general-fixed_obj}--\eqref{eq:general-fixed_argmin}
satisfying
\[
\ve z\in\proj_{\ve z}R_m \subseteq Q := \left\{\ve z \in \R^d:
    \begin{array}{lll}
        C\ve x^*+D\ve z\le\ve p, \ \ \ve z \ge \ve 0 , \ \ \ve c \cdot \ve x^*+\ve e\cdot \ve z = v^* \\
        \ve r \le B\ve z + \ve u < \ve r + \ve 1 \\
    \end{array}\right\}
\]
where $\ve 1$ is the all-ones vector in $\Z^m$.
As in the proof of Proposition~\ref{prop:denominator-bound} we note that every $\ve z \in Q$
gives rise to an optimal solution $(\ve x^*, \ve z)$.
We now show how to find one such optimal solution
whose binary encoding length is polynomial in the input size.
Let $k = 1 + \dim Q$.
Find $k$ affinely independent vertices $\ve z^1,\dots,\ve z^k$ of $\cl Q$.
This can be done in polynomial time since $\cl Q$ is a polytope
(see, e.g., \cite{Schrijver1986}, page 186).
By standard results in linear optimization, the
vertices $\ve z^1,\dots,\ve z^k$ can be written as rational numbers in
$\frac{1}{\Delta}\Z^d$ where $\Delta$ is
the least common multiple of the determinants of the bases of
the constraint
submatrix
associated with the basic feasible solutions $\ve z^1, \dots, \ve z^k$.
Note that $\Delta$ has encoding length polynomial in the input size.
It is possible that none of the $\ve z^j$
lie in $Q$, since $Q$ is not closed, but nonetheless
their barycenter
$\ve z^* = \sum_{j=1}^k \frac{1}{k}\; \ve z^j$
is guaranteed to be in $Q$ and also
in the set $\frac{1}{k\Delta}\Z^d$.
Thus it has a binary encoding length polynomial in the input size.

It follows that $(\ve x^*, \ve z^*)$ is
an optimal solution to \eqref{eq:general-fixed_obj}--\eqref{eq:general-fixed_argmin} and
can be found in polynomial time.
\end{proof}

We note, in particular, that this theorem provides a polynomial-time algorithm to decide if a BMILP
instance has an optimal solution.
To the authors' knowledge, this is the first algorithm to achieve this in the literature on discrete bilevel optimization. The existence of solutions in the case where $\ve z$ is continuous
has long been known not to be guaranteed (see for instance \cite{MooreBard1990})
but little is known in terms of necessary and sufficient conditions for existence.
Dempe and Richter \cite{Dempe2000-knapsack} explore the setting where the lower level problem is a knapsack problem, and they provide necessary and sufficient conditions for existence of solutions, but there have been no studies of this kind in a more general setting.
An efficient algorithm to decide the existence of solutions, such as the one in Theorem~\ref{thm:main-result-ES}, may be helpful in investigations of this kind.

A final question remains: what can we say about the BMILP problem \eqref{eq:general-fixed_obj}--\eqref{eq:general-fixed_argmin} when the infimum is not attained? In such
a case it may be desirable to find $\epsilon$-optimal solutions whose objective values approximate the infimum. The following result demonstrates how we can adjust our arguments to find $\epsilon$-optimal solutions.

\begin{corollary}\label{cor:epsilon-optima}
Let $\epsilon >0$ be given in binary encoding. Given a feasible instance of the BMILP problem \eqref{eq:general-fixed_obj}--\eqref{eq:general-fixed_argmin} with infimum value $v^*$, there exists an algorithm to find a feasible
solution $(\ve x, \ve z)$ that satisfies  $\ve c \cdot \ve x + \ve e \cdot \ve z \le (1+ \epsilon)v^*$. The algorithm runs in polynomial time when the follower's dimension $n$
is fixed.
\end{corollary}
\begin{proof}
By Theorem~\ref{thm:main-result-ES} part (ii) we can decide in polynomial time if the
infimum is attained. If so, by using part (iii) of the same theorem we can find an
optimal solution $(\ve x^*,\ve z^*)$ that satisfies $\ve c \cdot \ve x^* + \ve e \cdot \ve z^* = v^* \le (1+ \epsilon)v^*$.

If we find that the infimum is not attained we apply a similar procedure as
in the proof of Theorem~\ref{thm:main-result-ES} part (iii). By replacing the
condition $\ve c \cdot \ve x^*+\ve e\cdot \ve z = v^*$ with the condition
$\ve c \cdot \ve x + \ve e \cdot \ve z \le (1+ \epsilon)v^*$, the result follows immediately.
\end{proof}

It should be noted that this algorithm
achieves a stronger form
of efficiency than a fully polynomial time approximation scheme (\textsc{FPTAS}).
In an \textsc{FPTAS}, the error $\epsilon$ is given in unary encoding, whereas
in our algorithm $\epsilon$ enters as input in binary encoding,
since its value only plays a role in the constraint
$\ve c \cdot \ve x + \ve e \cdot \ve z \le (1+\epsilon)v^* + $.
In this sense, we have obtained
a ``better than fully polynomial time'' approximation scheme
to yield solutions
with objective value $\epsilon$-close to the infimum.

A distinguishing feature of our main result,
Theorem~\ref{thm:main-result-ES}, is that
the dimension $d$ of the leader's variable need not be fixed in order to assure a polynomial running time.
Note, however, that this result in no way conflicts with the fact that
bilevel \emph{linear} programming, in which both $\ve x$ and $\ve z$ are continuous
and the dimension of both
$\ve x$ and $\ve z$
are allowed to vary, is an NP-hard problem.
In our setting it is the leader's variable only which is allowed to be continuous and
have varying dimension, and an important ingredient in our proof is that
there are only finitely many
integer right-hand sides
$\beta_i(\ve z) = \lfloor B_i \ve z + \ve u \rfloor$ to consider
when the matrix $A$ and decision variables $\ve x$ are integer.

\section{An Algorithm for Bilevel Integer Linear Programming}\label{s:pure-integer}

We can also extend our results to the pure integer setting, where $\ve z$ is further restricted to be integer.
We refer to this problem as the bilevel integer linear programming (BILP) problem,
an instance of which is specified by the same data as BMILP and the same constraints
\eqref{eq:general-fixed_obj}--\eqref{eq:general-fixed_argmin} with the additional
integrality condition
\begin{equation}\label{eq:z-integer}
\ve z \in \Z^d.
\end{equation}
As above, we assume that the feasible set $\mathcal{F}'$
defined in~(\ref{eqn:pure_bilevel_feasible_optimistic})
is bounded.

\begin{theorem}\label{thm:BILP}
There exists an algorithm that solves
the BILP problem \eqref{eq:general-fixed_obj}--\eqref{eq:general-fixed_argmin}
\& \eqref{eq:z-integer} in the following sense: the algorithm
\begin{enumerate}[(i)]
\item decides if
a given
instance is feasible;
\item if it is feasible, finds an optimal solution.
\end{enumerate}
The algorithm runs in polynomial time when
the total dimension $d+n$ is fixed.
\end{theorem}
\begin{proof}
(i) Feasibility can be established in polynomial time by applying Lenstra's integer
programming algorithm \cite{Lenstra83} to decide if
$\{(\ve x, \ve z) \in P : A\ve x \le B\ve z + \ve u, \ve (\ve x, \ve z) \in \Z^{n+d}\} \neq \emptyset$.

(ii) Assume the instance is feasible.
The idea of how to determine an optimal solution is similar to,
but simpler than, the mixed integer case, as we use binary search to directly find
the lexicographically minimal optimal solution~$(\ve x^*, \ve z^*)$.
To begin, consider the decision
problem
$\BILP_\alpha$, which is identical to $\BMILP_\alpha$ with the
additional integrality condition~(\ref{eq:z-integer}).

Problem $\BILP_\alpha$ can be decided in polynomial time using an
algorithm similar to that
in Proposition~\ref{prop:decision-version}.
Indeed, if we simply redefine
\begin{eqnarray*}
Q_{\alpha} = \{ (\ve b, \ve x, \ve b) &:&
         \ve c\cdot\ve x + \ve e\cdot\ve b_{(3)}\le \alpha,
     \ \ C\ve x + D\ve b_{(3)} \le\ve p,
     \ \ \ve b_{(3)} \ge 0,             \\
     &&  A\ve x \le \ve b_{(2)},
     \ \ \ve b_{(1)} = \ve\psi\cdot\ve x - 1,
     \ \ \ve b_{(2)} = B\ve b_{(3)} + \ve u
  \},
\end{eqnarray*}
and let $p = n+d$,
the algorithm in the proof of Proposition~\ref{prop:decision-version} decides $\BILP_\alpha$.

It is then straightforward to find the optimal value $v^*$ of
the given BILP instance  using
binary search as in the proof of
Lemma~\ref{lemma:infimum-value}
and stopping
when the search interval has length less than one.

Given the optimal value $v^*$ we again use binary search to determine the lexicographically minimal optimal solution $(\ve x^*, \ve z^*)$. The procedure to find $\ve x^*$ is identical
to that described in the proof of Theorem~\ref{thm:main-result-ES}.
The procedure to find $\ve z^*$ is
also similar. We successively determine each component of $\ve z^*$ as follows:
assume that all components $z_k^*$ for $k = 1, \dots,j-1$ have already been
determined. Then $z_j^*$ is the optimal value of the BILP
\begin{eqnarray*}
    \min_{\ve z} && z_j  \\
    \st &&     C\ve x^* + D\ve z \le \ve p;\ \ \ve z \ge 0; \ \ \ve z \in \Z^d
            \ \ \ve c \cdot \ve x^* + \ve e \cdot \ve z = v^*; \\
        &&     z_k = z_k^*\;\text{ for }k = 1,\dots,j-1;        \\
    && \ve x^* \in \argmin_{\ve x'}\{\, \ve \psi \cdot \ve x':
            A \ve x' \leq B\ve z + \ve u,\ \ \ve x'\in\Z^n \,\}.
\end{eqnarray*}
Since this problem is feasible and its objective is integer valued and bounded, its
optimum value is attained, and can be determined as
in the previous paragraph.

The result then follows.
\end{proof}

We point out that an alternate algorithm satisfying the conditions of the previous theorem
can be obtained by applying a rational generating function technique similar to that
found in~\cite{Koppe2008-games}.
The algorithm is very similar to that in the Stackelberg--Nash setting of Section 6 in that paper.
The differences are:
(i) there is only a single follower as opposed to a group of followers
(hence any follower's optimal decision is trivially a ``Nash equilibrium'');
(ii) the payoffs are linear instead of piecewise linear; and
(iii) the set constraining the leader's actions
is more complicated and depends
on the action of the follower,
see \eqref{eq:general-fixed_P},
whereas in the Stackelberg--Nash setting
of~\cite{Koppe2008-games}
the leader's choice set is independent of the actions of the followers.
Nonetheless, these differences are minor and easily handled by
the rational generating function technique described in \cite{Koppe2008-games}.
Details are thus omitted.

How do these two algorithms to solve BILP -- the one based on parametric integer programming versus
the one based on rational generating function techniques -- compare in terms of efficiency?
It is accurate to say that the algorithm using parametric integer programming programming is more direct and
simpler.
Indeed, the rational generating function approach relies on a strong
result called the ``Projection Theorem" (Theorem~1.7 in \cite{Barvinok2003}) which,
although it has polynomial running time (in fixed dimension),
was until recently thought to be
virtually
unimplementable in practice.
Only recently has the Projection Theorem been implemented
in the library \texttt{barvinok} developed by Sven Verdoolaege \cite{barvinok}.
Moreover, the current implementation appeals to the results in \cite{Eisenbrand2008},
some of the very same results that we use
more directly
in the development of our first approach (see \cite{Koeppe2008-projection} for more details).
Future work would be required to improve implementations of the Projection Theorem
to yield a practical solution method for bilevel integer programs.

As for computational comparisons,
in the pure integer setting,
between the our two algorithms of Theorem~\ref{thm:BILP} herein and those
of Moore and Bard \cite{MooreBard1990} and DeNegre and Ralphs \cite{Denegre2008},
this is also a direction for future work.

\section{Conclusions}\label{s:conclusion}

In this paper we have detailed an algorithm which
solves the BMILP problem \eqref{eq:general-fixed_obj}--\eqref{eq:general-fixed_argmin}
in polynomial time when the dimension of the follower's integer decision variable is fixed.
It solves BMILP
in the following strong sense: under our boundedness and integral data assumptions,
it decides if the problem is feasible, and if so determines its infimum;
if the infimum is attained, it finds an optimal solution,
else it
finds an $\epsilon$-optimal solution for
any
$\epsilon > 0$ given as input in binary encoding.

This algorithm is based on
recent results in parametric integer programming due to Eisenbrand and Shmonin \cite{Eisenbrand2008}
combined with binary search to find the infimum objective value and an optimal or $\epsilon$-optimal solution.
%
A simplified version of our algorithm also solves the pure
integer problem BILP
\eqref{eq:general-fixed_obj}--\eqref{eq:general-fixed_argmin} \& \eqref{eq:z-integer}
in polynomial time when the total number of decision variables is fixed.

%
%
\bibliographystyle{plain}
\bibliography{references}

\end{document}